\documentclass[a4paper]{article}
\usepackage {diagrams, a4, theorem, amssymb, QED}
\title{Total Cofibres of Diagrams of Spectra}
\author{Thomas H\"uttemann}
\date{}

\newcommand{\bZ}{\mathbb{Z}}
\newcommand{\bN}{\mathbb{N}}
\newcommand{\iso}{\cong}
\newcommand{\pliso}{\iso_{PL}}
\newcommand{\id}{\mathrm{id}}

\newcommand{\sset}{\mathrm{sSet}}
\newcommand{\Fun}{\mathrm{Fun}\,}
\newcommand{\ie}{{\it i.\thinspace e.}}
\newcommand{\hocolim}{\mathrm{hocolim}}
\newcommand{\holim}{\mathrm{holim}}
\newcommand{\C}{\mathcal{C}}
\newcommand{\D}{\mathcal{D}}
\newcommand{\Sp}{\mathfrak{Sp}}
\newcommand{\Tsp}{\mathfrak{Tsp}}
\newcommand{\dash}{\hbox{-}\,}
\newcommand{\op}{\mathrm{op}}

\newarrow{Cof}>--->
\newarrow{Fib}----{>>}

\newtheorem{theorem}{Theorem}

\newtheorem{lemma}[theorem]{Lemma}

\theorembodyfont{\upshape}
\newtheorem{definition}[theorem]{Definition}
\newtheorem{remark}[theorem]{Remark}

\begin{document} 

\maketitle
\centerline {\it Universit\"at G\"ottingen, Fakult\"at f\"ur Mathematik}
\centerline {\it Mathematisches Institut, Bunsenstr.~3--5}
\centerline {\it D--37073 G\"ottingen, Germany}
\centerline {huette@uni-math.gwdg.de}
\centerline {http://www.uni-math.gwdg.de/huette/}
\medskip
\noindent{Published by New York Journal of Mathematics {\bf 11} (2005), pp.~333--343}\\
\noindent{Available online at: http://nyjm.albany.edu:8000/j/2005/11-16.html}

\begin{abstract}  
  If $Y$ is a diagram of spectra indexed by an arbitrary poset $\C$ together
  with a specified sub-poset $\D$, we define the {\it total cofibre} $\Gamma
  (Y)$ of~$Y$ as
  $$\mathrm{cofibre}(\mathrm{hocolim}_\D (Y) \rTo \mathrm{hocolim}_\C (Y)) \ .$$
  We construct a comparison map $\hat\Gamma_Y \colon
  \holim_\C Y \rTo \hom (Z, \hat\Gamma (Y))$ to a mapping spectrum of
  a fibrant replacement of
  $\Gamma (Y)$ where $Z$ is a simplicial set obtained from~$\C$ and~$\D$, and
  characterise those poset pairs $\D \subset \C$ for which $\hat\Gamma_Y$ is a
  stable equivalence. The characterisation is given in terms of stable
  cohomotopy of spaces related to~$Z$. For example, if $\C$ is a finite
  polytopal complex with $|\C| \iso B^m$ a ball with boundary sphere $|\D|$,
  then $|Z|\pliso S^m$, and $\hat\Gamma(Y)$ and $\holim_\C (Y)$ agree up to
  $m$-fold looping and up to stable equivalence.  As an application of the
  general result we give a spectral sequence for $\pi_*(\Gamma(Y))$ with
  $E_2$-term involving higher derived inverse limits of $\pi_* (Y)$,
  generalising earlier constructions for space-valued diagrams indexed by the
  face lattice of a polytope.
\end{abstract} 

\tableofcontents

\section{Introduction}

Let $\D \subseteq \C$ be a pair of posets, considered as categories (with
arrows starting at the smaller element). The motivating example will be
the poset $\C(P)$ of non-empty faces of a polytope~$P$, and its sub-poset
$\D(P)$ of proper faces. More generally, we will consider polytopal
complexes $\D \subset \C$ with underlying spaces $PL$-homeomorphic to $S^{m-1}
\subset B^m$.

\begin{definition}
  \label{def:gamma}
  Let $X \colon \C \rTo \sset_*,\ F \mapsto X^F$ be a diagram of pointed
  simplicial sets.  The {\it total cofibre\/} $\Gamma (X)$ of the diagram $X$
  is the (strict) cofibre of the map $\hocolim_\D (X) \rTo \hocolim_\C (X)$.
\end{definition}

The homotopy colimit is the coend $\hocolim_\C (X) = N(\dash \downarrow \C)_+
\otimes_\C X$ formed with respect to the smash product of pointed simplicial
sets. In a similar way, the total cofibre functor can be defined for diagrams
of pointed topological spaces (or, more generally, for diagrams with values in
a pointed simplicial model category).

For diagrams with values in pointed topological spaces, indexed by $\C =
\C(P)$ and $\D = \D(P)$, there is a convergent first quadrant spectral
sequence
$$E^2_{p,q} = \lim{}^{n-p} \tilde H_q (X; \bZ) \,\Longrightarrow\, \tilde
H_{p+q} (\Gamma (X); \bZ)$$
where $n = \dim P$, cf.~\cite[Theorem~2.23]{H-Finiteness} (note that the
definition of~$\Gamma$ given there is homotopy equivalent to the one given
here since the natural map 
$$\hocolim_\D Y \rTo \hocolim_\C Y$$
is a cofibration, hence cofibre and homotopy cofibre are homotopy equivalent). The
appearance of higher derived {\it inverse\/} limits in the $E^2$-term might be
unexpected and results from an argument using \textsc{\v Cech} cohomology.
This paper offers a different and more satisfying explanation: We
prove that, stably, total cofibres and homotopy limits agree. More precisely, we
construct a comparison map from homotopy limits to a loop space of the total
cofibre (Definition~\ref{def:comp-map}), and prove that under certain
conditions on $\D$ and $\C$ this map is a weak equivalence in stable homotopy
(\ie, for diagrams of spectra instead of simplicial sets;
Theorem~\ref{thm:main}).  Then the usual \textsc{Bousfield}-\textsc{Kan}
spectral sequence for homotopy limits can be used to construct a spectral
sequence of the type above (Theorem~\ref{thm:bkss}). The conditions on $\C$
and $\D$ are phrased in terms of stable cohomotopy and
completely characterise those poset pairs for which total cofibres and
homotopy limits agree.

The total cofibre functor has been used by the author to examine categories of
space-valued quasi-coherent sheaves on projective toric varieties
\cite{H-Finiteness, H-ksplit}. It replaces the global sections functor and its
higher deriveds, the sheaf cohomology functors, in algebraic geometry. The
spectral sequence mentioned above provides a means to relate the topological
construction to the algebraic context \cite[Theorem~3.8]{H-Finiteness} since
computation of $\lim^p$ amounts to taking $p$th sheaf cohomology. Note also
that computing sheaf cohomology can be thought of as a homotopy limit process
in a category of chain complexes, so the results of this paper show that in
fact the total cofibre functor is completely analogous to sheaf cohomology.

\section{Spectra and stable cohomotopy}
\label{sec:conventions}

We denote the category of \textsc{Bousfield}-\textsc{Friedlander} spectra
\cite{BF:Gamma} by~$\Sp$. Thus a spectrum $X$ is a sequence of pointed
simplicial sets $X_0, X_1, \ldots$ and structure maps $\Sigma X_n \rTo
X_{n+1}$.  The category $\Sp$ has two simplicial closed model structures, the
{\it level structure\/} with levelwise weak equivalences and levelwise
fibrations, and the {\it stable structure\/} with stable equivalences and with
cofibrations as in the level structure \cite[\S2]{BF:Gamma}.

Both model structures of~$\Sp$ are cofibrantly generated
\cite[\S11.1]{Hh-Modelcats} with generating cofibrations the maps
$\alpha_{j,k} \colon F_j (\partial \Delta^k_+) \rTo F_j (\Delta^k_+)$ for $j,k
\in \bN$. Here $F_j (A)$ is the free spectrum generated by $A \in \sset_*$ in
degree~$j$; in other words, $F_j(A)$ is the suspension spectrum of~$A$ shifted
$j$ times. The cofibre of $\alpha_{j,k}$ is the spectrum $\tilde K_{j,k} :=
F_j(\Delta^k/\partial\Delta^k)$. We let $\tilde K_{j,k} \rCof^\simeq K_{j,k}$
denote a stably fibrant replacement.

For $W \in \Sp$ and $A \in \sset_*$ we denote by $\hom_{\sset_*} (A,W)$ the
spectrum which has the mapping space $\hom_{\sset_*} (A, W_n)$ in level $n$;
since the hom functor commutes with looping, we can use the structure maps
$W_n \rTo \Omega W_{n+1}$ to make $\hom_{\sset_*} (A,W)$ into a spectrum.

\begin{definition}
  \label{def:cohomotopy}
  For $n \in \bZ$ the $n$th stable cohomotopy group of $A \in \sset_*$ is
  defined as $\pi_s^n (A):= \pi_m \hom_{\sset_*} (A, K_{j,k})$ where $m \in
  \bZ$ and $j,k \in \bN$ are numbers satisfying $k-j-m = n$.
\end{definition}

The definition of $\pi_s^n (A)$ does not depend on the choices of $m$, $j$ and
$k$: There are isomorphisms $\pi_m \hom_{\sset_*} (A, K_{j,k}) \iso \pi_{m-1}
\hom_{\sset_*} (S^1,\, \hom_{\sset_*} (A, K_{j,k}))\allowbreak\iso
\pi_{m-1} \hom_{\sset_*} (A,\, \hom_{\sset_*} (S^1, K_{j,k}))$, and the two
spectra $\hom_{\sset_*} (S^1, K_{j,k})$ and $K_{j+1, k}$ are weakly
equivalent. Similarly, $\hom_{\sset_*} (S^1, K_{j,k}) \simeq K_{j,k-1}$.
Similar arguments work for suspension instead of looping.

For a $\C$-diagram $X$ in $\sset_*$ recall that $\holim_\C (X) = \hom_\C (N(\C
\downarrow \dash)_+, X)$, the space of natural transformations $N(\C
\downarrow \dash)_+ \rTo X$. If $Y$ is a diagram of spectra we can define
$\holim_\C (Y) = \hom_\C (N(\C \downarrow \dash)_+, Y)$ by forming the
homotopy limit in each level; since $\holim$ commutes with looping, we can use
the structure maps $Y_n \rTo \Omega Y_{n+1}$ to make the result into a
spectrum.

\section{Comparing homotopy limits and total cofibres}

Let $Y \colon \C \rTo \Sp$ be a diagram of spectra.  Since the
functor $\Gamma$ (Definition~\ref{def:gamma}) commutes with suspensions, we
can apply $\Gamma$ in each level of~$Y$ and obtain a spectrum $\Gamma (Y)$.
Thus we can and will consider $\Gamma$ as a functor from diagrams of spectra
to spectra.

\begin{lemma}
\label{lem:Gamma_hococart}
  The functor $\Gamma$ commutes with homotopy colimits. If $Y \rTo Z$ is a
  natural transformation of diagrams of spectra consisting of stable
  equivalences, then $\Gamma (Y) \rTo \Gamma(Z)$ is a stable equivalence.
\end{lemma}

\begin{proof}
  The first assertion is immediate from the definition. For the second, note
  that homotopy colimits preserve stable equivalences
  \cite[Lemma~5.18]{Th-etale}. The claim then follows by comparing the long exact
  sequence of stable homotopy groups associated to the injection $\iota \colon \hocolim_\D
  Y \rTo \hocolim_\C Y$ and the corresponding sequence for~$Z$ in place
  of~$Y$. (Note that the homotopy cofibre of~$\iota$ is levelwise weakly
  equivalent to the strict cofibre, and similarly for~$Z$ instead of~$Y$.)
\end{proof}

\begin{lemma}\label{lem:cont}
  The functor $\Gamma$ is continuous (\ie, induces a map of $\hom$-spaces).
\end{lemma}

\begin{proof}
  Let $Y,Z \colon \C \rTo \Sp$ be two diagrams.  An $n$-simplex on the mapping
  space $\hom_\C (Y,Z)$ is a natural transformation of the form $Y \wedge
  \Delta^n_+ \rTo Z $.  Since $\Gamma$ commutes with the functor $\dash \wedge
  \Delta^n_+$, application of the total cofibre functor yields a map $\Gamma
  (Y) \wedge \Delta^n_+ \rTo \Gamma (Z)$ which is a simplex in the mapping
  space $\hom_\Sp (\Gamma(Y), \Gamma(Z))$.
\end{proof}

\begin{definition}\label{def:comp-map}
  Let $Y \colon \C \rTo \Sp$ be a diagram of spectra. By Lemma~\ref{lem:cont}
  the functor $\Gamma$ induces a natural map of spectra
  $$\Gamma_Y \colon \holim_{\C} (Y) = \hom_{\C} \big( N(\C \downarrow
  \dash)_+, Y \big) \rTo \hom_{\sset_*} \big( \Gamma (N(\C \downarrow
  \dash)_+), \Gamma (Y) \big) \ .$$
%%% NEW
  Let $j_Y \colon \Gamma (Y) \rCof^\simeq \hat\Gamma (Y)$ denote a functorial
  level fibrant replacement; we denote by $\hat\Gamma_Y$ the composition
  $$\hat\Gamma_Y = (j_Y)_* \circ \Gamma_Y \colon \holim_{\C} (Y) \rTo
  \hom_{\sset_*} \big( \Gamma (N(\C \downarrow
  \dash)_+), \hat\Gamma (Y) \big) \ .$$
%%% ENDNEW
\end{definition}

We need another bit of notation. Suppose that $\C$ is a (possibly infinite)
poset, considered as a category; for any $F \in \C$ we write $\C^F := \{ G \in
\C \,|\, F \not\leq G \}$. If $\C$ happens to be a polytopal or simplicial
complex, then $\C^F$ is $\C$ with the open star of~$F$ removed (see
\cite[\S1.1]{H-ksplit} for more on this terminology).

\bigbreak

Suppose $\D$ is an order ideal in the poset $\C$ (\ie, $\D \subseteq \C$, and
for all $x,y \in \C$, if $x \leq y \in \D$ then $x \in \D$). We will consider
the following conditions on the pair $(\D,\C)$:
\begin{list}{\rm(P\theenumi)}{\usecounter{enumi}}
\item The space $N(\C)/N(\C^F)$ has trivial cohomotopy for all $F \in \D$,
  \ie, for all $F \in \D$ and all $m \in \bZ$ we have $\pi_s^m (N(\C)/N(\C^F)) = 0$.
\item The map $\beta \colon N(\C)/N(\D) \rTo N(\C)/N(\C^F)$ is a stable
  cohomotopy equivalence for all $F \in \C \setminus \D$; equivalently, the
  homotopy cofibre of~$\beta$ has trivial stable cohomotopy for all $F \in \C
  \setminus \D$.
\end{list}

\begin{theorem}\label{thm:main}
  Suppose $\D$ is an order ideal in the poset $\C$. The following statements
  are equivalent:
  \begin{enumerate}
  \item Given any diagram $Y \colon \C \rTo \Sp,\ F \mapsto Y^F$ of spectra
    with $Y^F$ levelwise \textsc{Kan} for all $F \in \C$, the comparison map
    $${\hat\Gamma_Y} \colon \holim_{\C} (Y) \rTo \hom_{\sset_*} \big( \Gamma (N(\C
    \downarrow \dash)_+), \hat\Gamma (Y) \big)$$
    is a stable weak equivalence of spectra.
  \item The pair $(\D,\C)$ satisfies the conditions (P1) and~(P2) specified above.
  \end{enumerate}
\end{theorem}

\begin{lemma}\label{lem:special}
  Suppose that the posets $\D$ and $\C$ satisfy the conditions (P1) and~(P2).
  Fix $j,k \in \bN$.  Given $F \in \C$ we define a diagram of spectra
  $$Q := \C (F,\,\dash)_+ \wedge K_{j,k} \colon \C \rTo \Sp, \quad G \mapsto
  \C(F,\,G)_+ \wedge K_{j,k} $$
  where $K_{j,k}$ is as defined in \S\ref{sec:conventions}.  Then $\Gamma_Q$
  and $\hat\Gamma_Q$ are weak equivalences.
\end{lemma}

\begin{proof}
{\it Case~1: $F \in \D$.\/} Recall that $\Gamma (Q)$ is the (strict)
cofibre of the map
$$ \hocolim_{\D} (Q) \rTo^\kappa \hocolim_{\C} (Q) \ .$$
Since $Q = \C (F,\,\dash)_+ \wedge K_{j,k}$, we have isomorphism
$$ \hocolim_{\D} Q \iso N (F \downarrow \D)_+ \wedge K_{j,k}
\quad\hbox{and}\quad \hocolim_{\C} Q \iso 
N(F \downarrow \C)_+ \wedge K_{j,k} \ , $$
the map~$\kappa$ being induced by the inclusion $\D \subset \C$. Both $F
\downarrow \D$ and $F \downarrow \C$ have the initial object $F$, hence their
nerves are simplicially contractible, with contraction induced by the functor
sending every object to $F$. Thus $\Gamma (Q)$ is simplicially contractible as
well. Since mapping spaces respect simplicial homotopies we conclude
$$\hom_{\sset_*} (\Gamma(N(\C\downarrow \dash)_+),\, \Gamma(Q)) \simeq * \ .$$
%%% NEW
Now $\Gamma (Q) \simeq *$ implies $\hat\Gamma (Q) \simeq *$, hence
$\hom_{\sset_*} (\Gamma(N(\C\downarrow \dash)_+),\, \hat\Gamma(Q)) \simeq *$
since $\hat\Gamma (Q)$ is level fibrant.
%%% ENDNEW

On the other hand, by definition $ \holim_\C (Q) = \hom_\C \left(N(\C
  \downarrow \dash)_+, Q \right)$. Since $Q$ is the constant diagram with
value $K_{j,k}$ on $F \downarrow \C$, and since $Q (G) = *$ if $G \in \C^F$
we have an isomorphism
$$ \holim_{\C} (Q) \iso \hom_{\sset_*} (N(\C)/N(\C^F), K_{j,k}) \ .$$
The target spectrum has homotopy groups the stable cohomotopy of the space
$N(\C)/N(\C^F)$ (Definition~\ref{def:cohomotopy}).  By hypothesis~(P1) this
means that both sides of the comparison map
$$\Gamma_Q \colon \holim_{\C} (Q) \rTo \hom_{\sset_*} (\Gamma(N(\C\downarrow \dash)_+),\,
\Gamma(Q))$$
are weakly contractible, 
%%% NEW
and the same is true for $\hat\Gamma_Q$.
%%% ENDNEW
Hence $\Gamma_Q$ and $\hat\Gamma_Q$ are weak equivalences as claimed.

{\it Case~2: $F \in \C \setminus \D$.\/} Then $Q$ is trivial on~$\D$ and hence
 $\hocolim_\D (Q) = *$. The natural map~$\alpha$ from homotopy colimit to colimit
 $$\alpha \colon \hocolim_\C (Q) \iso \Gamma (Q) \iso N(F \downarrow \C)_+
 \wedge K_{j,k} \rTo \mathrm{colim}_\C Q \iso S^0 \wedge K_{j,k} \iso K_{j,k}$$
 is induced by
 sending all non-basepoints in $N(F \downarrow \C)_+$ to the non-basepoint
 in~$S^0$. Note that $N(F \downarrow \C)$ is simplicially contractible, hence
 $\alpha$~is a simplicial homotopy equivalence.

Arguing as in Case~1 we see that $ \holim_{\C} Q \iso \hom_{\sset_*}
(N(\C)/N(\C^F),\, K_{j,k}) $. Moreover, since $\D$ is an order ideal in~$\C$ we have
$\D \subseteq \C^F$. Thus we obtain the following commutative diagram:
\begin{diagram}
  \holim_\C Q = \hom_\C \left(N(\C \downarrow \dash)_+, Q\right) & \rTo[l>=3em]^{\Gamma_Q} &
  \hom_{\sset_*} \left(\Gamma (N(\C \downarrow \dash)_+), \Gamma (Q)\right) \\
  \dTo<{\iso} && \\
  \hom_{\sset_*} \left(N(\C)/N(\C^F), K_{j,k}\right) && \dTo>\simeq<{\alpha_*} \\
  \dTo<f && \\
  \hom_{\sset_*} \left(N(\C)/N(\D), K_{j,k}\right) & \rTo_g & \hom_{\sset_*} \left(\Gamma (N(\C
  \downarrow \dash)_+), K_{j,k}\right)
\end{diagram}
Here $\alpha_*$ is a weak equivalence since $\alpha$ is a simplicial homotopy
equivalence. The map~$f$ is a weak equivalence by~(P2).
The map $g$ is induced by the
natural weak equivalences from homotopy colimits to colimits
$$\hocolim_\D N(\D \downarrow \dash)_+ \rTo^\simeq N\D_+ \quad\hbox{and}\quad
\hocolim_\C N(\C \downarrow \dash)_+ \rTo^\simeq N\C_+ \ ,\eqno{(*)}$$
hence is a weak equivalence as well; recall that since $K_{j,k}$ is fibrant, the
$\hom$-functor preserves weak equivalences.  (To see why the maps~$(*)$ are
weak equivalences, one can argue as follows.  For any category $\mathcal{E}$
we have (using {\it unpointed topological spaces\/} for ease of notation, and
writing $\hocolim^\prime$ for unpointed homotopy colimits)
$$\hocolim^\prime_\mathcal{E} |N(\mathcal{E} \downarrow \dash)| = |N (\dash
\downarrow \mathcal{E})| \otimes_\mathcal{E} |N (\mathcal{E} \downarrow
\dash)| \rTo^{(a)} \hocolim^\prime_{\mathcal{E}^\op} |N (\dash \downarrow \mathcal{E})| $$
$$\qquad\qquad \rTo^\simeq \hocolim^\prime_{\mathcal{E}^\op} |*| \rTo^{(a)}
|*| \otimes_\mathcal{E} |N(\mathcal{E} \downarrow \dash)| =
\mathrm{colim}_\mathcal{E} |N(\mathcal{E} \downarrow \dash)| = |N\mathcal{E}|
\ ,$$
the arrow marked with ``$\simeq$'' being a weak equivalence by homotopy
invariance of homotopy colimits. The two arrows marked $(a)$ use the {\it
  anti-simplicial\/} identification $N (\C \downarrow \dash) \iso N (\dash
\downarrow \C^\op)$, inducing homeomorphisms after geometric realisation. The
composite map is induced by the canonical map from homotopy colimits to
colimits, as can be checked by a direct calculation.)

Since four of the maps in the diagram are weak equivalences, the fifth map
$\Gamma_Q$ is necessarily a weak equivalence as well.

%%% NEW
Finally, since $K_{j_k}$ is level fibrant we can choose a map $\ell \colon
\hat\Gamma(Q) \rTo K_{j,k}$ with $\ell \circ j_Q = \alpha$
(cf.~Definition~\ref{def:comp-map}). Then $\ell$ is a level equivalence
between level fibrant spectra, hence $\ell_*$ is a weak equivalence of hom spectra.
From $\ell_* \circ (j_Q)_* = \alpha_*$ we conclude that $(j_Q)_*$ is a weak
equivalence, hence so is $\hat\Gamma_Q = (j_Q)_* \circ \Gamma_Q$.
%%% ENDNEW
\end{proof}

\bigbreak

\noindent {\bf Proof of Theorem \ref{thm:main}.}
The category $\Fun (\C, \Sp)$ has a model structure where a morphism $Y \rTo
X$ is a weak equivalence (fibration) if and only if for all $F \in \C$ the map
$Y^F \rTo X^F$ is a levelwise weak equivalence (levelwise fibration). We call
the resulting model structure on $\Fun (\C, \Sp)$ the {\it $\ell$-structure\/}
and will speak of $\ell$-equivalences, $\ell$-fibrations, {\it etc}. We will
use the symbol $\simeq$ to denote $\ell$-equivalences.  Similarly, the
category $\Fun (\C, \Sp)$ has a model structure where a morphism $Y \rTo X$ is
a weak equivalence (fibration) if and only if for all $F \in \C$ the map $Y^F
\rTo X^F$ is a stable weak equivalence (stable fibration).  We call the
resulting model structure on $\Fun (\C, \Sp)$ the {\it $s$-structure\/} and
will speak of $s$-equivalences, $s$-fibrations, {\it etc}. We will use the
symbol $\sim$ to denote $s$-equivalences. Cofibrations are the same in both
cases.

Suppose that $\D$ and~$\C$ satisfy conditions (P1) and~(P2).  Let $Y \colon \C
\rTo \Sp$ be $\ell$-fibrant, and let $Y^\prime \rTo^\simeq Y$ be an
$\ell$-cofibrant replacement of~$Y$. Then $Y^\prime$ is $\ell$-fibrant.
Homotopy colimits map $\ell$-equivalences to level equivalences since homotopy
colimits of diagrams of simplicial sets are weakly homotopy invariant.
Homotopy limits map $\ell$-equivalences of $\ell$-fibrant objects to
$\ell$-equivalences since homotopy limits of diagrams of \textsc{Kan} sets are
weakly homotopy invariant. From the commutative square below we infer that we
may assume $Y$ to be $\ell$-cofibrant and $\ell$-fibrant.

\begin{diagram}[textflow]
  \holim (Y^\prime) & \rTo[l>=4em]^{\hat\Gamma_{Y^\prime}} & \hom_{\sset_*}
  (\Gamma (N(\C \downarrow \dash)_+),
  \hat\Gamma (Y^\prime)) \\
  \dTo<\simeq && \dTo>\simeq \\
  \holim (Y) & \rTo^{{\hat\Gamma_Y}} & \hom_{\sset_*} (\Gamma (N(\C \downarrow
  \dash)_+), \hat\Gamma (Y)) \\
\end{diagram}

Since the level structure of $\Sp$ is cofibrantly generated, so is the
$\ell$-structure of $\Fun (\C, \Sp)$; for each $F \in \C$ and each generating
cofibration $i \colon A \rTo B$ of $\Sp$ we have one generating cofibration
$\C (F,\dash)_+ \wedge A \rTo^{\id \wedge i} \C (F, \dash)_+ \wedge B$ for the
$\ell$-structure. Since $Y$ is $\ell$-cofibrant, it is a retract of an
$\ell$-cellular object. By naturality of $\hat\Gamma_Y$ it is thus enough to prove
that $\hat\Gamma_Y$ is a stable equivalence for $Y$ an $\ell$-cellular object
(which might not be $\ell$-fibrant).

Since $Y$ is cellular we may filter $Y$ by a (possibly transfinite) sequence
$$ \ldots \rTo Y_i \rTo^{\alpha_i} Y_{i+1} \rTo^{\alpha_{i+1}} \ldots
\eqno{(\dagger)}$$
starting with the trivial diagram $Y_0 = *$, where $\alpha_i$ is a pushout of
a generating cofibration. (This is a special case of the small object
argument, cf.~\cite[Corollary~11.2.2]{Hh-Modelcats}.)  In particular
$\alpha_i$ is a cofibration and $Y_i$ is cofibrant. (If $\C$ is countable the
sequence $(\dagger)$ is indexed by the natural numbers. In general, it will be
indexed by a regular cardinal $\lambda$, and at each limit ordinal $\kappa <
\lambda$ the natural map $\mathrm{colim}_{\mu < \kappa} Y_\mu \rTo Y_\kappa$
is an isomorphism.)

\begin{figure}[htbp]
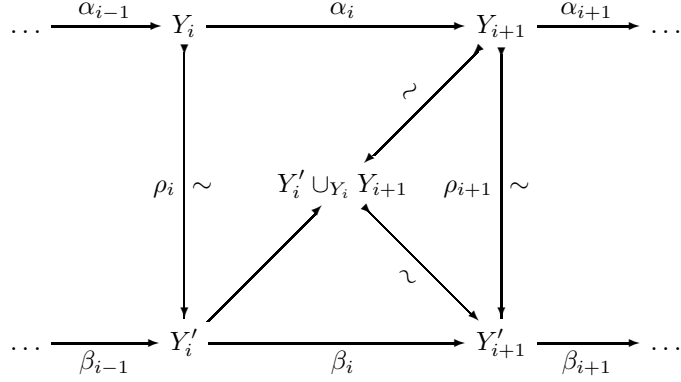

  \centering
  \begin{diagram}
    \ldots & \rTo^{\alpha_{i-1}} &Y_i&& \rTo^{\alpha_i}&&Y_{i+1}&
    \rTo^{\alpha_{i+1}}& \ldots \\
    &&&&& \ldCof<\sim \\
    && \dCof<{\rho_i}>\sim && Y_i^\prime \cup_{Y_i} Y_{i+1} && \dCof<{\rho_{i+1}}>\sim \\
    &&& \ruTo && \rdCof<\sim \\
    \ldots & \rTo_{\beta_{i-1}} &Y_i^\prime&& \rTo_{\beta_i} &&Y_{i+1}^\prime&
    \rTo_{\beta_{i+1}}& \ldots \\
  \end{diagram}
  \caption{Replacement in the \textsc{Reedy} model structure}
  \label{fig:reedy}
\end{figure}

We can replace this sequence by an $s$-equivalent one consisting of objects
which are $s$-cofibrant and $s$-fibrant. This can be done by fibrant
replacement in the \textsc{Reedy\/} model structure of sequences of diagrams,
formed with respect to the $s$-structure \cite[\S15.3]{Hh-Modelcats}. More
explicitly, set $Y_0^\prime := Y_0 = *$ and $\rho_0 := \id_{Y_0}$. Given $\rho_i
\colon Y_i \rCof^\sim Y_i^\prime$ we define $\rho_{i+1} \colon Y_{i+1}
\rCof^\sim Y_{i+1}^\prime$ as the pushout of $\rho_i$ along~$\alpha_i$
composed with an $s$-fibrant replacement $Y_i^\prime \cup_{Y_i} Y_{i+1}
\rCof^\sim Y_{i+1}^\prime$; see Fig.~\ref{fig:reedy}.
(For a limit ordinal $\kappa$ we let $\rho_\kappa$
denote the composition of $\mathrm{colim}_{\mu < \kappa} \rho_\mu$ composed
with an $s$-fibrant replacement $\mathrm{colim}_{\mu < \kappa} Y_\mu^\prime
\rCof^\sim Y_\kappa^\prime$.) Note that by construction $Y_i^\prime$ is fibrant and
cofibrant in both the $s$-structure and the $\ell$-structure, and the natural
maps $\beta_i \colon Y_i^\prime \rTo Y_{i+1}^\prime$ are cofibrations.

\noindent We then certainly have a chain of weak equivalences
$$\mathrm{hocofibre}\, (\beta_i) \lTo^\sim \mathrm{hocofibre}\, (\alpha_i)
\rTo^\simeq \mathrm{cofibre}\, (\alpha_i) \ .$$
Since $\alpha_i$ is a generating
cofibration, its cofibre is of the form $\C(F_i ,\dash)_+ \wedge \tilde
K_{j,k}$ for some element $F_i \in \C$ and numbers $j,k \in \bN$ (see
\S\ref{sec:conventions} for the definition of $\tilde K_{j,k}$ and $K_{j,k}$).
The stably fibrant replacement $\tilde K_{j,k} \rCof^\simeq K_{j,k}$ induces
an $s$-equivalence
$$\mathrm{cofibre}\, (\alpha_i) \rTo^\sim \C(F_i, \dash)_+ \wedge K_{j,k} =: Q_i$$
to an $s$-fibrant diagram~$Q_i$. Since $\mathrm{hocofibre}\, (\beta_i)$ is
$s$-cofibrant, we can find a single $s$-equivalence
$$\omega_i \colon \mathrm{hocofibre}\, (\beta_i) \rTo^\simeq Q_i $$
representing the zig-zag constructed above.  By composing this with the
canonical map $Y_{i+1}^\prime \rTo \mathrm{hocofibre}\, (\beta_i)$ we get a
sequence
$$Y_i^\prime \rTo^{\beta_i} Y_{i+1}^\prime \rTo Q_i \ .$$
By construction this is a cofibration sequence in the sense that the composite
map is weakly null homotopic (in the homotopy category the composite factors
as $Y_i^\prime \rTo Y_i \rTo \mathrm{cofibre}\, \alpha_i \rTo Q_i$ which is
the trivial map), and the map $\mathrm{hocofibre}\, (\beta_i) \rTo^{\omega_i}
Q_i$ is an $s$-equivalence.  Since $Y_i^\prime$ and $Q_i$ are $s$-cofibrant
and $s$-fibrant, the composite map is actually simplicially null homotopic and
factors over the simplicial reduced cone $ CY_i^\prime := Y_i^\prime \wedge
\Delta^1$ of~$Y_i^\prime$ (here $\Delta^1$ has basepoint~$1$). Let $E_i$
denote an $s$-fibrant replacement of the cone $CY_i^\prime$. By the lifting
axioms, we can extend the map $CY_i^\prime \rTo Q_i$ to~$E_i$, so we get a
commutative diagram of $s$-fibrant and $s$-cofibrant objects
\begin{diagram}
Y_i^\prime & \rTo^{\beta_i} & Y_{i+1}^\prime \\
\dTo && \dTo \\
E_i & \rTo & Q_i \\
\end{diagram}
which is homotopy cocartesian, \ie, for each $F \in \C$ the canonical map
$$ \mathrm{hocolim} \left( E_i^F \lTo Y_i^F \rTo Y_{i+1}^F \right) \rTo Q_i^F $$
is a stable equivalence. By the properties of spectra, the square is then homotopy
cartesian as well, \ie, for each $F \in \C$ the canonical map
$$ Y_i^F \rTo \mathrm{holim} \left( E_i^F \rTo Q_i^F \lTo Y_{i+1}^F \right) $$
is a stable equivalence.

We now use that if $\hat\Gamma_{Y_i^\prime}$ and $\hat\Gamma_{Q_i}$ are stable
equivalences, then so is $\hat\Gamma_{Y_{i+1}^\prime}$. This follows from the fact
that $\holim$ preserves homotopy cartesian squares, and that $\hat\Gamma$
preserves homotopy cocartesian squares (the latter follows immediately from
Lemma~\ref{lem:Gamma_hococart}, the former follows from a formally dual
argument).
%%% NEW
Moreover, the hom space functor preserves homotopy cartesian squares of level
fibrant spectra.
%%% ENDNEW
Hence the comparison maps $\hat\Gamma_?$ assemble to a map of square
diagrams of spectra
%%% NEW
(writing $Z = \Gamma (N(\C \downarrow \dash)_+)$ for
compact notation)
%%% ENDNEW
$${
  \begin{diagram}
    \holim (Y_i^\prime) & \rTo & \holim (Y_{i+1}^\prime) \\ \dTo && \dTo \\
    \holim (E_i) & \rTo & \holim (Q_i) \\
  \end{diagram} \rTo[l>=4em]^{\hat\Gamma_?}
  \begin{diagram}
    \hom_{\sset_*} (Z, \hat\Gamma
    (Y_i^\prime)) & \rTo &
    \hom_{\sset_*} (Z, \hat\Gamma
    (Y_{i+1}^\prime)) \\
    \dTo && \dTo \\
    \hom_{\sset_*} (Z, \hat\Gamma (E_i)) &
    \rTo &
    \hom_{\sset_*} (Z, \hat\Gamma (Q_i)) \\
  \end{diagram}
}$$
where source and target are homotopy cartesian. The map is a weak equivalences
on three corners; this is true for $\hat\Gamma_{Y_i^\prime}$ and
$\hat\Gamma_{Q_i}$ by assumption, and since $E_i$ is simplicially
contractible, we have $\holim (E_i) \simeq * \simeq \hat\Gamma (E_i)$. Hence
the map of diagrams a weak equivalence everywhere.

Consequently, since $\hat\Gamma_{Y_0}$ is a weak equivalence , it is
sufficient to check that $\hat\Gamma_{Q_i}$~is a weak equivalence. But this
follows from Lemma~\ref{lem:special}; the diagram~$Q$ defined there is now
called $Q_i$.

\medbreak

To prove the other implication, suppose first that condition~(P1) is not
satisfied. Then we can find $F \in \D$ and $j,k \in \bN$ such that
$\hom_{\sset_*} (N(\C)/N(\C^F), K_{j,k})$ is not contractible, and the the
proof of Lemma~\ref{lem:special} shows that $Q = \C(F,\,\dash)_+ \wedge
K_{j,k}$ is a diagram of spectra such that $\hat\Gamma_Q$ is not an
equivalence of spectra. A similar argument works if (P2) is not satisfied.
\qed

\begin{remark}
\label{rem:strongcond}
  Conditions~(P1) and~(P2) are implied by the following stronger conditions:
    \begin{list}{(P\theenumi')}{\usecounter{enumi}}
    \item For every $F \in \D$ the inclusion $N\C^F \rTo N\C$ is a weak
      equivalence.
    \item For every $F \in \C \setminus \D$ the inclusion $N \D \rTo N \C^F$
      is a weak equivalence.
    \end{list}
    For (P2') implies that the the map $N(\C)/N(D) \rTo N(\C)/N(\C^F)$ is a
    weak equivalence, and (P1') implies $N(\C)/N(\C^F) \simeq *$. Conditions
    (P1) and~(P2) now follow directly from weak homotopy invariance of the
    mapping spectra $\hom_{\sset_*} (\dash,\,K_{j,k})$.
\end{remark}

\section{The Bousfield-Kan spectral sequence}

\subsection*{Diagrams defined on finite posets}

\begin{remark}
  For finite posets, conditions~(P1) and~(P2) can be simplified by applying
  the following fact to the geometric realisation of $\mathrm{hocofibre}
  (\beta)$ and the simplicial set~$N\C/N\C^F$: {\it A finite pointed $CW$
    complex $K$ has trivial stable cohomotopy if and only if $\Sigma^2 K$ is
    contractible.}  Indeed, $\Sigma^2 K \simeq *$ clearly implies $\pi_s^n (K)
  = 0$ for all~$n$.  Conversely, if $K$ has trivial stable cohomotopy, then
  its \textsc{Spanier}-\textsc{Whitehead} dual $DK$ has trivial stable
  homotopy, hence is stably contractible. Thus the double dual $DDK$ is stably
  contractible, hence has trivial reduced homology. Since $K$ is finite there
  is a homotopy equivalence $K \simeq DDK$. So $K$ has trivial homology, hence
  so has $\Sigma^2 K$. Consequently $\Sigma^2 K$, being simply connected, is
  contractible.
\end{remark}

\begin{theorem}\label{thm:bkss}
  Let $\D \subseteq \C$ be a pair of finite posets satisfying the conditions~(P1)
  and~(P2), and let $Y \colon \C \rTo \Sp$ be a diagram with $Y^F$ levelwise
  \textsc{Kan} for all $F \in \C$. There is a strongly
  convergent right-halfplane spectral sequence
  $$ E_2^{p,q} = \lim{}^p (\pi_q Y) \,\Longrightarrow\, \pi_{q-p}
  \hom_{\sset_*}\Big(\Gamma (N(\C\downarrow\dash)_+), \hat\Gamma (Y)\Big) \ ;$$
  the $r$-differential is a map $d_r \colon E_r^{p,q} \rTo E_r^{p+r,
  q+r-1}$.
\end{theorem}

\begin{proof}
  If $X \in \Sp$ is levelwise \textsc{Kan}, we can define a new spectrum~$QX$
  by setting $(QX)_n := \mathrm{colim}_{k\rightarrow\infty} \Omega^k X_{n+k}$.
  The spectrum $QX$ depends functorially on~$X$; it is stably fibrant, and the
  natural map $X \rTo QX$ is a stable equivalence.
  
  By forming $QY^F$ for all $F \in \C$ we obtain a diagram $QY$ consisting
  entirely of stably fibrant spectra, and an $s$-equivalence $Y \rTo QY$. Now
  since $\C$ is a finite poset, all the over-categories $\C \downarrow F$ have
  finite nerve. Since $\hom_{\sset_*} (K,\, \dash)$ commutes with filtered
  colimits for any finite simplicial set~$K$, the functor $Q$ commutes with
  $\holim$ and the composite map
  $$ \holim_\C (Y) \rTo Q \holim_\C (Y) \iso \holim_C (QY) \eqno{(*)}$$
  is a stable equivalence. Now by \cite[Proposition~5.13]{Th-etale} we have a
  spectral sequence
  $$ E_2^{p,q} = \lim{}^p (\pi_q Y) \,\Longrightarrow\, \pi_{q-p} (\holim_\C
  (QY)) $$
  which converges strongly since $N\C$ is finite-dimensional. Using the stable
  equivalence~$(*)$ and Theorem~\ref{thm:main} we obtain a natural isomorphism
  $$\pi_{q-p} (\holim_\C (QY)) \iso \pi_{q-p} \hom_{\sset_*}\Big( \Gamma
  (N(\C\downarrow\dash)_+), \hat\Gamma(Y) \Big) \ . $$
\end{proof}

\subsection*{Topological spectra}

In all we did so far we could replace the category $\Sp$ of
\textsc{Bousfield}-\textsc{Friedlander} spectra with the category $\Tsp$ of
topological spectra (sequences of pointed topological spaces $(X_i)_{i \in
  \bN}$ and structure maps $\Sigma X_i \rTo X_{i+1}$). Indeed, the category
$\Tsp$ has a (topological) model structure similar to the one of~$\Sp$
\cite[\S2.5]{BF:Gamma}; a set of generating cofibrations consists of the maps
$$|\alpha_{j,k}| \colon F_j(|\partial \Delta^k_+|) \rTo F_j(|\Delta^k_+|)$$
where $F_j(T)$ is the suspension spectrum of~$T$ shifted by~$j$. We can also
define the functor $\Gamma$ and the map
$$\Gamma_Y \colon \holim_{\C} (Y) = \hom_{\C} \big( |N(\C \downarrow
\dash)_+|, Y \big) \rTo \hom \big( \Gamma (|N(\C \downarrow \dash)_+|), \Gamma
(Y) \big)$$
in this context. The statements of Theorem~\ref{thm:main},
Lemma~\ref{lem:special} and Theorem~\ref{thm:bkss}
%%% NEW
with $\hat\Gamma$ replaced by $\Gamma$
%%% ENDNEW
remain valid {\it mutatis
  mutandis}; in fact their proofs could be simplified slightly since all
topological spectra are automatically levelwise fibrant. This latter
observation implies that {\it given $K \in \sset_*$ and $Y \in \Tsp$, if $|K|$ is
homeomorphic to an $m$-sphere there is a natural isomorphism $\pi_p \hom (|K|,
Y) \iso \pi_{p+m} Y$\/}. Here $\hom (|K|, Y)$ is the topological spectrum
which has the topological mapping space $\hom (|K|, Y_n)$ in level~$n$.

\subsection*{Diagrams defined on $PL$ balls}

A {\it polytopal complex\/} $\C$ is a collection of non-empty polytopes in
$\mathbb{R}^n$ which is closed under taking non-empty faces, such that the
intersection of $F,G \in \C$ is a (possibly empty) face of~$F$ and~$G$. For
any subset $\mathcal{K} \subseteq \C$ we write $|\mathcal{K}|:= \bigcup_{K \in
  \mathcal{K}} K$ for the underlying space of~$\mathcal{K}$.

Let $\C$ be (the underlying poset of) a finite polytopal complex with $|\C|$ a
$PL$ ball of dimension~$m$, and let $\D$ be the subcomplex with $|\D| =
\partial |\C|$. Then $N\C$ and $N\D$ are the barycentric subdivisions of $\C$
and $\D$, respectively. For $F \in \C\setminus\D$ the poset $\C^F$ is obtained
from $\C$ by removing the open star $st_\C (F)$ of~$F$ in~$\C$, \ie, removing
all polytopes containing~$F$. Thus $|N\C^F|$ is the $m$-ball $|N\C|$ with an
open $m$-ball in the interior removed, hence condition~(P2') is satisfied.
Similarly, condition~(P1') holds since for $F \in \D$ the open star $st_\C(F)$
of $F$ in~$\C$ is an $m$-ball, and the open star of $F$ in~$\D$ is an
$(m-1)$-ball with underlying space $|st_\C (F)| \cap |\D|$.

By Remark~\ref{rem:strongcond}, Theorem~\ref{thm:main} applies to diagrams $\C
\rTo \Tsp$. Since there is a $PL$ homeomorphism of pairs
$$\left( |\hocolim_\C (N(\C\downarrow\dash)_+)|,\ |\hocolim_\D
  (N(\D\downarrow\dash)_+)| \right) \,\iso\, (|N\C_+|,\ |N\D_+|)\ ,$$
and since $ |\C|/|\D| \iso S^m $, this means that total cofibres and homotopy
limits agree up to $m$-fold looping and up to stable weak equivalences:

\begin{theorem}
  Let $\C$ be a finite polytopal complex with $|\C|$ a $PL$ ball of
  dimension~$m$, and let $\D$ be the subcomplex with $|\D| = \partial |\C|$.
  If $Y \colon \C \rTo \Tsp$ is a diagram of topological spectra, the natural
  map
  $$\Gamma_Y \colon \holim_\C (Y) \rTo \hom \big( \Gamma (|N(\C \downarrow
  \dash)|_+), \Gamma (Y) \big)$$
  is a stable equivalence of topological spectra.
  Moreover, Theorem~\ref{thm:bkss} yields a strongly convergent right
  half-plane spectral sequence
  $$E_2^{p,q} = \lim{}^p (\pi_q Y) \,\Longrightarrow\, \pi_{m+q-p} (\Gamma
  (Y)) \ . $$
  \qed
\end{theorem}

If $X \colon \C \rTo \mathrm{Top}_*$ is a diagram of pointed topological
spaces, we can consider the diagram of topological spectra $X \wedge H\bZ$
where $H\bZ$ denotes the \textsc{Eilenberg}-\textsc{Maclane} spectrum
of~$\bZ$. Then $\pi_q (X \wedge H\bZ) \iso \tilde H_q (X; \bZ)$ (reduced
integral homology), and we obtain the following theorem generalising
\cite[Theorem~2.23]{H-Finiteness} (where the result is stated only for the
special case $\C = \C(P)$ and $\D = \D(P)$, and a different indexing
convention is used for the spectral sequence):

\begin{theorem}
  Let $\C$ be a finite polytopal complex with $|\C|$ a $PL$ ball of
  dimension~$m$, and let $\D$ be the subcomplex with $|\D| = \partial |\C|$.
  For any diagram $X \colon \C \rTo \mathrm{Top}_*$ there is a strongly convergent
  spectral sequence
  $$E_2^{p,q} = \lim{}^p (\tilde H_q (X; \bZ)) \,\Longrightarrow\, \tilde
  H_{m+q-p} (\Gamma (X); \bZ)\ . $$
  \qed
\end{theorem}

\subsection*{Acknowledgements}
The author is indebted to \textsc{S.~Schwede} and \textsc{M.~Weiss} for
helpful comments. All diagrams typeset with Paul Taylor's ``diagrams'' macro
package for \LaTeX.

%%%%%%%%%%%%%%%%%%%%%%%%%%%%%%%%%%%%%%%%%%%%%%%%%%%%%%%%%%%%%%%%%%

\ifx\undefined\bysame
\newcommand{\bysame}{\leavevmode\hbox to3em{\hrulefill}\,}
\fi

\end{document}